\documentclass[12pt,reqno,intlimits]{amsart}
\usepackage[cp1251]{inputenc}
\usepackage[T2A]{fontenc}
\usepackage{amsfonts,amsmath,amsxtra,amsthm,amssymb,latexsym}

\theoremstyle{plain}
\newtheorem {theorem} {Theorem}[section]

\theoremstyle{definition}
\newtheorem {definition}[theorem]{Definition}
\newtheorem {remark}[theorem]{Remark}

\begin{document}

\title{Elimination of Hamilton--Jacobi equation in extreme variational problems}
\author{Igor Orlov}
 \maketitle
\centerline{Taurida National V. Vernadsky University, Simferopol,
Ukraine} \centerline{E-mail address: old@crimea.edu}

\begin{abstract}
It is shown that extreme problem for one--dimensional
Euler--Lagrange variational functional in ${C^1[a;b]}$ under the
strengthened Legendre condition can be solved without using
Hamilton--Jacobi equation. In this case, exactly one of the two
possible cases requires a restriction to a length of $[a;b]$,
defined only by the form of integrand. The result is extended to
the case of compact extremum in ${H^1[a;b]}$.

\end{abstract}

\section*{\textbf{Introduction}}

The classical scheme of the research to a local extremum for a
one--dimensional Euler--Lagrange variational functional
$$
\Phi(y)=\int\limits_a^b f(x,y,y')dx \mapsto \mbox{extr}\quad(y\in
C^1[a;b])
$$
at an extremal point $y$ assumes~\cite{oiv_1},~\cite{oiv_2}
checking the strengthened Legendre condition
${f_{y'y'}(x,y,y')\neq 0}$ and the Jacobi condition ${U(x)\neq 0}$
${(a<x\leq b)}$ for the Hamilton--Jacobi equation:
$$
-\frac{d}{dx}\biggl[f_{y'y'}(x,y,y')
U'\biggr]+\biggl[-\frac{d}{dx}\bigl(f_{yy'}(x,y,y')\bigr)+f_{y^2}(x,y,y')
\biggr]U=0\quad (U(a)=0,\,U'(a)=1).
$$

The second step is the most laborious, it requires to solve a
complicated enough equation with a view to receive, really, a very
small information about behavior of the solution $U(x)$.

Moreover, the initial conditions ${U(a)=0,\,U'(a)=1}$, have as a
consequence automatical fulfilment of the Jacobi condition near
$a$. The question is only --- how much length has a suitable
interval?

The aim of the present work is to show that the interval
satisfying the Jacobi condition can be chosen depending only on
the form of the integrand $f$ and not depending on a concrete
extremal. More precisely, the main result
(Theorem~\ref{Th_oiv_1.1},~\ref{Th_oiv_3.1}) distinguishes two
cases depending on the range of the coefficients in the
Hamilton--Jacobi equation. For the first case, an extremum is
guaranteed without any restriction to a length of $[a;b]$, for the
second one, such a restriction is presented. The result above
remains valid under the passage to the case of the research to a
compact extremum in Sobolev space ${H^1[a;b]}$.

The first part of the work deals with elimination of the
Hamilton--Jacobi equation in case of zero extremal in
${C^1[a;b]}$. The second part contains a quadratic estimate of
tending $\Phi$ to a minimal value via norm of $y$ in ${H^1[a;b]}$.
The third part determines a general form of $\Phi$ under the
conditions of a local minimum and Legendre condition at zero. The
fourth part contains a passage to the case of an arbitrary
$C^2$--smooth extremal in ${C^1[a;b]}$ and the last, fifth part
contains a passage to the case of a compact minimum in
${H^1[a;b]}$.

\section{\textbf{Elimination of Jacobi condition: case of zero extremal}}

Let's consider a classical Euler--Lagrange variational functional
\begin{equation}\label{oiv_1}
\Phi(y)=\int\limits_a^b f(x,y,y')dx\quad (y\in
C^1[a;b],\,y(a)=y(b)=0,\,f\in C^2,\, f_{yz}\in C^1).
\end{equation}

We are going to show that, under fulfilment of Euler--Lagrange
variational equation and the strengthened Legendre condition at
zero, the functional~\eqref{oiv_1} \textit{always} attains a
strong local extremum at zero. However, in addition, two different
possible cases defined by form of the integrand $f$, arise: one of
the cases assumes a restriction to a length of ${[a;b]}$, at the
second case any restriction is absent.

So, let's divide the integrand ${f(x,y,z)}$ into two terms:
$$
f_1(x,y,z)=f(x,y,z)-f(x,0,0)-\left[f_y(x,0,0)\cdot
y+f_z(x,0,0)\cdot
z\right]-
$$
$$
-\frac{1}{2}\left[f_{y^2}(x,0,0)\cdot
y^2+2f_{yz}(x,0,0)\cdot yz+\lambda\cdot f_{z^2}(x,0,0)\cdot
z^2\right];\quad (0<\lambda<1)
$$
$$
f_2(x,y,z)=f(x,y,z)-f_1(x,y,z)=f(x,0,0)+\left[f_y(x,0,0)\cdot
y+f_z(x,0,0)\cdot
z\right]+
$$
$$
+\frac{1}{2}\left[f_{y^2}(x,0,0)\cdot
y^2+2f_{yz}(x,0,0)\cdot yz+\lambda\cdot f_{z^2}(x,0,0)\cdot
z^2\right].
$$

\quad

Let's set, respectively,
$$
\Phi_i(y)=\int\limits_a^b f_i(x,y,y')dx\quad
(i=1,2);\quad\Phi(y)=\Phi_1(y)+\Phi_2(y).
$$

1) Let's investigate $\Phi_1$ for a local extremum (minimum, for
definiteness) at zero with the help of Euler--Lagrange, Legendre
and Jacobi conditions.

(i) \textit{The Euler--Lagrange equation.} Because
$$
\bigl(f_{1,y}(x,y,z)=f_y(x,y,z)-f_y(x,0,0)-f_{y^2}(x,0,0)\cdot
y-f_{yz}(x,0,0)\cdot z\bigr)\Rightarrow
$$
$$
\quad\quad\quad\quad\Rightarrow\bigl(f_{1,y}(x,0,0)=0\bigr);
$$
$$
\bigl(f_{1,z}(x,y,z)=f_z(x,y,z)-f_z(x,0,0)-f_{yz}(x,0,0)\cdot
y-\lambda\cdot f_{z^2}(x,0,0)\cdot z\bigr)\Rightarrow
$$
$$
\quad\quad\quad\quad\Rightarrow\bigl(f_{1,z}(x,0,0)=0\bigr);
$$
then the Euler--Lagrange equation for $\Phi_1$ at zero
$$
f_{1,y}(x,0,0)-\frac{d}{dx}\bigl[f_{1,z}(x,0,0)\bigr]=0
$$
holds automatically, i.e. $y_0(x)\equiv 0$ is an extremal of the
functional $\Phi_1$.

(ii) \textit{The strengthened Legendre condition.} Because
$$
\bigl(f_{1,z^2}(x,y,z)=f_{z^2}(x,y,z)-\lambda\cdot
f_{z^2}(x,0,0)\bigr)\Rightarrow
\bigl(f_{1,z^2}(x,0,0)=(1-\lambda)\cdot f_{z^2}(x,0,0)\bigr),
$$
then, under the additional requirement
\begin{equation}\label{oiv_2}
p(x):=f_{z^2}(x,0,0)>0,\quad (a\leq x\leq b)
\end{equation}
the strengthened Legendre condition for a strong minimum at zero
holds.

(iii) \textit{The Hamilton--Jacobi equation and the Jacobi
condition.} Because
$$
\bigl(f_{1,yz}(x,y,z)=f_{yz}(x,y,z)-f_{yz}(x,0,0)\bigr)\Rightarrow
\bigl(f_{1,yz}(x,0,0)=0\bigr);
$$
$$
\bigl(f_{1,y^2}(x,y,z)=f_{y^2}(x,y,z)-f_{y^2}(x,0,0)\bigr)\Rightarrow
\bigl(f_{1,y^2}(x,0,0)=0\bigr);
$$
then the Hamilton--Jacobi equation for $\Phi_1$ at zero takes form
of
$$
-\frac{d}{dx}\biggl[(1-\lambda)\cdot
f_{z^2}(x,0,0)U'\biggr]+\left[-\frac{d}{dx}
\bigl(f_{1,yz}(x,0,0)\bigr)+f_{1,y^2}(x,0,0)\right]U=
$$
$$
=-\frac{d}{dx}\biggl[(1-\lambda)p(x)U'\biggr]=0\quad\,\,(U(a)=0,\,\,U'(a)=1).
$$

Hence, in view of condition~\eqref{oiv_2}, the required result
$$
\biggl(U(x)=p(a)\cdot
\int\limits_a^x\frac{dt}{p(t)}\biggr)\Rightarrow\biggl(U(x)\neq 0
\,\,\mbox{for}\,\, a<x\leq b\biggr)
$$
holds, i.e. the strengthened Jacobi condition at zero for a strong
minimum of $\Phi_1$ takes place. Thus, under the
condition~\eqref{oiv_2}, $\Phi_1$ attains a strong local minimum
at zero.

2) Let's investigate now $\Phi_2$ for a local extremum at zero
immediately. Note at first that $\Phi_2(0)=\Phi(0)$.

(i) Suppose that the Euler--Lagrange equation for $\Phi$ at zero
\begin{equation}\label{oiv_3}
f_y(x,0,0)-f_{xz}(x,0,0)=0\quad (a\leq x\leq b)
\end{equation}
holds. Then integrating by parts gives us
$$
\Phi_2(y)=\int\limits_a^b
f(x,0,0)dx+\int\limits_a^b\left[f_y(x,0,0)\cdot y+f_z(x,0,0)\cdot
y'\right]dx+
$$
$$
+\int\limits_a^b\left[\frac{1}{2}f_{y^2}(x,0,0)\cdot
y^2+f_{yz}(x,0,0)\cdot
yy'\right]dx+\frac{\lambda}{2}\cdot\int\limits_a^b\cdot
f_{z^2}(x,0,0)\cdot y'^2 dx=
$$
$$
=\Phi_2(0)+\left[\int\limits_a^b\left(f_y-f_{xz}\right)(x,0,0)dx+
f_z(x,0,0)\cdot y\biggl|\biggr._a^b\right]+
$$
$$
+\left[\frac{1}{2}\int\limits_a^b\left(f_{y^2}-f_{xyz}\right)(x,0,0)\cdot
y^2 dx+ \frac{1}{2}f_{yz}(x,0,0)\cdot
y^2\biggl|\biggr._a^b\right]+ \frac{\lambda}{2}\int\limits_a^b
p(x)\cdot y'^2 dx.
$$
From here, denoting by
$$
q(x):=\left(f_{y^2}-f_{xyz}\right)(x,0,0),
$$
it follows
\begin{equation}\label{oiv_4}
 \Phi_2(y)=\Phi_2(0)+\frac{1}{2}\int\limits_a^b\left[\lambda\cdot p(x)\cdot
 y'^2+q(x)\cdot y^2\right]dx.
\end{equation}

(ii) Denote by
\begin{equation}\label{oiv_5}
p:=\min\limits_{a\leq x\leq b}p(x)>0,\quad\quad
q:=\min\limits_{a\leq x\leq b}q(x),
\end{equation}
and consider at first the case of ${q\geq 0}$. Then
$$
\lambda p(x) y'^2+q(x) y^2\geq \lambda p \cdot y'^2+q\cdot
y^2>0\quad \mbox{as}\quad y'\neq 0,
$$
whence, in view of~\eqref{oiv_4}, the inequality
$$
\Phi_2(y)>\Phi_2(0)\quad \mbox{as}\quad y(x)\neq 0
$$
follows. Thus, in this case $\Phi_2$ attains a strong absolute
minimum at zero. Hence, in view of one was proved in i.1), $\Phi$
attains a strong local minimum at zero (without any restriction to
a length of ${[a;b]}$).

(iii) Let's consider now the case of ${q<0}$. Then, using
Friederichs inequality (see, e.g.,~\cite{oiv_3}, Ch.~18), it
follows
$$
\Phi_2(y)-\Phi_2(0)=\frac{1}{2}\int\limits_a^b\left[\lambda\cdot
p(x)\cdot y'^2+q(x)\cdot
y^2\right]dx\geq\frac{1}{2}\int\limits_a^b\left[\lambda\cdot
p\cdot y'^2-|q|\cdot y^2\right]dx\geq
$$
\begin{equation}\label{oiv_6}
\geq \frac{1}{2}\int\limits_a^b\left[\lambda\cdot p\cdot
y'^2-\frac{16(b-a)^2}{\pi^2}|q|\cdot
y'^2\right]dx=\frac{1}{2}\left(\lambda\cdot
p-\frac{16(b-a)^2}{\pi^2}|q|\right)\cdot\int\limits_a^b y'^2dx.
\end{equation}

Let's require that the coefficient in front of the last integral
in~\eqref{oiv_6} will be strictly positive:
\begin{equation}\label{oiv_7}
\left(\lambda\cdot
p-\frac{16(b-a)^2}{\pi^2}|q|>0\right)\Leftrightarrow\left(b-a<\frac{\pi}{4}\sqrt{\frac{\lambda
p}{|q|}}\right).
\end{equation}
It follows from~\eqref{oiv_6} and~\eqref{oiv_7} that
${\Phi_2(y)>\Phi_2(0)\quad \mbox{as}\quad y\neq 0}$, i.e. $\Phi_2$
attains a strong absolute minimum at zero and hence, by virtue of
one was proved in i.1), $\Phi$ attains a strong local minimum at
zero under the restriction~\eqref{oiv_7} to a length of ${[a;b]}$.

Finally, passing to the limits in~\eqref{oiv_7} as
${\lambda\rightarrow 1-0}$, the last statement can be extended to
the case of the estimate of a length of ${[a;b]}$ not depending on
$\lambda$:
$$
b-a<\frac{\pi}{4}\sqrt{\frac{p}{|q|}}.
$$

So, it is proved the following
\begin{theorem}\label{Th_oiv_1.1}
Let the variational functional~\eqref{oiv_1} satisfies at zero the
Euler--Lagrange equation~\eqref{oiv_3} under the conditions
${y(a)=y(b)=0}$. Then, under the notation of~\eqref{oiv_5},
\begin{itemize}
\item[1)] for ${p>0}$, ${q\geq 0}$, $\Phi(y)$ attains a strong
local minimum at zero (without any restriction to a length of
${[a;b]}$);

\item[2)]for ${p>0}$, ${q< 0}$, under the restriction to a
length of ${[a;b]}$:
\begin{equation}\label{oiv_8}
b-a<\frac{\pi}{4}\sqrt{\frac{p}{|q|}},
\end{equation}
$\Phi(y)$ attains a strong local minimum at zero as well.
\end{itemize}

\end{theorem}

\section{\textbf{Quadratic estimation from below of tending $\Phi$ to minimum at zero}}

It's easy to see that the estimate~\eqref{oiv_8} at
Theorem~\ref{Th_oiv_1.1} is not optimal. For example, a
generalized harmonic oscillator
$$
\Phi(y)=\int\limits_0^T(py'^2-qy^2)dx\quad\,\,(p>0,\,\,q>0)
$$
on zero extremal reduces to the Hamilton--Jacobi equation
$$
pU''+qU=0\quad\,\,(U(0)=0,\,\,U'(0)=1)
$$
having the solution
$$
U(x)=\sqrt{\frac{p}{q}}\sin\sqrt{\frac{q}{p}}x\,,
$$
satisfying Jacobi condition ${U(x)\neq 0}$ as ${o<x<T}$ for
${T<\pi\sqrt{\frac{p}{q}}}$\,.

At the same time, the estimate~\eqref{oiv_8} for given case leads
to inequality ${T<\frac{\pi}{4}\sqrt{\frac{p}{q}}}$\,. However, as
it's easily can be seen, an advantage of the estimate
~\eqref{oiv_8} consists of possibility to get a useful quadratic
estimate from below for tending $\Phi(y)$ to the minimal value by
means of norm of $y$ in the Sobolev space ${H^1[a;b]}$.

1) First, let's consider a case of ${p>0}$, ${q>0}$. The
equality~\eqref{oiv_4} implies
$$
\Phi_2(y)-\Phi_2(0)\geq \frac{1}{2}\min(p,q)\cdot
\int\limits_a^b(y'^2+y^2)dx=\frac{1}{2}\min(p,q)\cdot
\|y\|^2_{H^1[a;b]}\,\, .
$$
Since ${\Phi(y)-\Phi(0)\geq \Phi_2(y)-\Phi_2(0)}$ in a small
enough zero neighborhood, then given a zero neighborhood the
inequality
$$
\Phi(y)-\Phi(0)\geq \frac{1}{2}\min(p,q)\cdot \|y\|^2_{H^1[a;b]}
$$
holds true.

2) Let's pass to the case of ${p>0}$, ${q<0}$. The
inequality~\eqref{oiv_6} leads to the estimate
$$
\Phi_2(y)-\Phi_2(0)\geq \frac{1}{2}\left[
p-\frac{16(b-a)^2}{\pi^2}|q|\right]\cdot \int\limits_a^by'^2dx\,\,
.
$$
Since the Friederichs inequality implies
\begin{equation}\label{oiv_9}
\int\limits_a^by'^2dx\geq \frac{\pi^2}{\pi^2+16(b-a)^2} \cdot
\|y\|^2_{H^1[a;b]}\,\, ,
\end{equation}
then by combining of the last two inequalities for a small enough
neighborhood of zero, under the conditions of
inequality~\eqref{oiv_6}, we get
$$
\Phi(y)-\Phi(0)\geq
\frac{\pi^2p-16(b-a)^2|q|}{2(\pi^2+16(b-a)^2)}\cdot
\|y\|^2_{H^1[a;b]}\, .
$$

3) Note that the estimate~\eqref{oiv_9} can be applied as well in
the case of ${p>0}$, ${q\geq 0}$, whence the inequality
$$
\Phi(y)-\Phi(0)\geq \frac{\pi^2p}{2(\pi^2+16(b-a)^2)}\cdot
\|y\|^2_{H^1[a;b]}
$$
follows. So, it is proved the following

\begin{theorem}\label{Th_oiv_2.1}
Under the conditions and notation of Theorem~\ref{Th_oiv_1.1}, the
following statements are valid:
\begin{itemize}
\item[1)] in the case of ${p>0}$, ${q> 0}$, in small enough zero
neighborhood in $C^1[a;b]$ the estimate
$$
\Phi(y)-\Phi(0)\geq \frac{1}{2}\min(p,q)\cdot \|y\|^2_{H^1[a;b]}
$$
holds;

\item[2)] in the case of ${p>0}$, ${q\geq 0}$, in small enough
zero neighborhood in $C^1[a;b]$ the estimate
$$
\Phi(y)-\Phi(0)\geq \frac{\pi^2p}{2(\pi^2+16(b-a)^2)}\cdot
\|y\|^2_{H^1[a;b]}
$$
holds;

\item[3)] in the case of ${p>0}$, ${q<0}$, in small enough zero
neighborhood in $C^1[a;b]$, under the condition of
estimate~\eqref{oiv_8}, the estimate
$$
\Phi(y)-\Phi(0)\geq
\frac{\pi^2p-16(b-a)^2|q|}{2(\pi^2+16(b-a)^2)}\cdot
\|y\|^2_{H^1[a;b]}
$$
holds.
\end{itemize}
\end{theorem}

\section{\textbf{Application: inverse extreme problem for variational functional}}

Let's set up a problem: to find a general form of the variational
functional~\eqref{oiv_1} possessing local minimum at zero under
the strengthened Legendre condition.

1) Let's shall find an integrand $f$ of the
functional~\eqref{oiv_1} in the form of
\begin{equation}\label{oiv_10}
f(x,y,z)=P(x,y)+Q(x,y)\cdot z+\frac{1}{2}R(x,y,z)\cdot z^2\, .
\end{equation}
Then
$$
P(x,y)=f(x,y,0),\quad Q(x,y)=f_z(x,y,0),\quad
R(x,y,0)=f_{z^2}(x,y,0).
$$
Under this notation, the Euler--Lagrange equation on zero
extremal~\eqref{oiv_3} takes form of
\begin{equation}\label{oiv_11}
(Q_x-P_y)(x,0)=0\quad (a\leq x\leq b);
\end{equation}
the strengthened Legendre condition on zero extremal~\eqref{oiv_2}
takes form of
\begin{equation}\label{oiv_12}
R(x,0,0)=:p(x)>0\quad (a\leq x\leq b).
\end{equation}

2) Let's choose an arbitrary ${P(x,y)\in C^2}$. Then a general
form of $Q$ follows from~\eqref{oiv_11}:
$$
\biggl(Q_x(x,0)=P_y(x,0)\biggr)\Rightarrow
\biggl(Q(x,0)=C+\int\limits_a^x P_y(t,0)dt\biggr)\Rightarrow
$$
$$
\Rightarrow \biggl(Q(x,y)=C+\int\limits_a^x
P_y(t,0)dt+\widetilde{Q}(x,y),\,\,\,\mbox{where}\,\,\,\widetilde{Q}(x,0)=0\biggr)\Rightarrow
$$
\begin{equation}\label{oiv_13}
\Rightarrow \biggl(Q(x,y)=C+\int\limits_a^x
P_y(t,0)dt+[q(x,y)-q(x,0)]\biggr),
\end{equation}
here ${C\in \mathbb{R}}$ and ${q(x,y)\in C^2}$ can be chosen
arbitrarily.

3) A general form of $R$ easily follows from the
condition~\eqref{oiv_12}:
\begin{equation}\label{oiv_14}
\biggl(R(x,0,0)=p(x)>0\biggr)\Rightarrow
\biggl(R(x,y,z)=p(x)+[\rho(x,y,z)-\rho(x,0,0)]\biggr),
\end{equation}
where ${p(x)>0}$, ${p\in C^2}$; ${\rho(x,y,z)\in C^2}$ can be
chosen arbitrarily.

4) A general form of the integrand $f$ follows now
from~\eqref{oiv_10},~\eqref{oiv_13} and~\eqref{oiv_14}:
$$
f(x,y,z)=P(x,y)+\biggl(C+\int\limits_a^x
P_y(t,0)dt+[q(x,y)-q(x,0)]\biggr)\cdot z+
$$
\begin{equation}\label{oiv_15}
+\frac{1}{2}\biggl(p(x)+[\rho(x,y,z)-\rho(x,0,0)]\biggr)\cdot
z^2\, ,
\end{equation}
where ${C\in\mathbb{R}}$; ${q,\,p\in C^2}$ ${(p>0)}$ can be chosen
arbitrarily. So, it is proved the following

\begin{theorem}\label{Th_oiv_3.1}
Let, under the conditions of Theorem~\ref{Th_oiv_1.1}, the
functional~\eqref{oiv_1} attains a local minimum at zero under the
strengthened Legendre condition. Then the integrand $f$ takes form
of~\eqref{oiv_15}.
\end{theorem}

\begin{remark}\label{Rem_oiv_3.1}
As it follows from Theorem~\ref{Th_oiv_3.1}, a general form of the
variational functional~\eqref{oiv_1} taking a local minimum at
zero under the strengthened Legendre condition is
$$
\Phi(y)=\int\limits_a^b \biggl(P(x,y)+\biggl[\int\limits_a^x
P_y(t,0)dt+q(x,y)-q(x,0)\biggr]\cdot y'+
$$
\begin{equation}\label{oiv_16}
+\frac{1}{2}\biggl[p(x)+\rho(x,y,y')-\rho(x,0,0)\biggr]\cdot
y'^2\biggr)dx\, ,
\end{equation}
where $P$, $q$, ${p>0}$, $\rho$ are the arbitrary functions from
$C^2$.

Thus, under the strengthened Legendre condition, the inverse
extreme variational problem at zero is solved: all the functionals
of type~\eqref{oiv_1} taking a local minimum at zero are
described.
\end{remark}

\section{\textbf{Case of arbitrary $C^2$--smooth extremal in $C^1[a;b]$}}

Let's fix an arbitrary $C^2$--smooth function $y_0(x)$, ${a\leq
x\leq b}$, and consider a question on elimination of Jacobi
condition for the local minimum of the variational
functional~\eqref{oiv_1} at the point $y_0(\cdot)$ under the
boundary conditions ${y(a)=y_0(a)}$, ${y(b)=y_0(b)}$.

To pass to the considered above (i.1) case of zero extremal, it
suffices to consider an auxiliary variational functional:
$$
\widetilde{\Phi}(y)=\Phi(y+y_0)=\int\limits_a^b
f(x,y+y_0(x),y'+y'_0(x))dx=:\int\limits_a^b
\widetilde{f}(x,y,y')dx
$$
$$
(y(a)=y(b)=0).
$$
In this connection the condition ${y_0(\cdot)\in C^2}$ guarantees
fulfilment of the condition from~\eqref{oiv_1} for the auxiliary
integrand $\widetilde{f}$ and permits to apply
Theorem~\ref{Th_oiv_1.1} to $\widetilde{\Phi}$. A not complicated
calculation shows that it is valid the following

\begin{theorem}\label{Th_oiv_4.1}
Let variational functional~\eqref{oiv_1} satisfies at a point
${y_0(\cdot)\in C^2[a;b]}$ Euler--Lagrange equation
\begin{equation}\label{oiv_17}
f_y(x,y_0,y'_0)-\frac{d}{dx}\bigl[f_z(x,y_0,y'_0)\bigr]=0\, .
\end{equation}
Denote by
$$
p:=\min\limits_{a\leq x\leq b}f_{z^2}(x,y_0(x),y'_0(x));
$$
$$
q:=\min\limits_{a\leq x\leq
b}\left[f_{y^2}(x,y_0(x),y'_0(x))-\frac{d}{dx}\bigl(f_{yz}(x,y_0(x),y'_0(x))\bigr)\right].
$$
Then, under the boundary conditions ${y(a)=y_0(a)}$,
${y(b)=y_0(b)}$,
\begin{itemize}
\item[1)] for ${p>0}$, ${q\geq 0}$, $\Phi(y)$ attains a strong
local minimum at $y_0(\cdot)$ (without any restriction to a length
of ${[a;b]}$);

\item[2)]for ${p>0}$, ${q<0}$, and under the restriction
\begin{equation}\label{oiv_18}
b-a<\frac{\pi}{4}\sqrt{\frac{p}{|q|}},
\end{equation}
to a length of ${[a;b]}$, $\Phi(y)$ attains a strong local minimum
at $y_0(\cdot)$ as well.
\end{itemize}
\end{theorem}

Analogously, applying Theorem~\ref{Th_oiv_2.1} to
$\widetilde{\Phi}$ leads a general quadratic estimate for tending
$\Phi$ to a local minimum at $y_0$.

\begin{theorem}\label{Th_oiv_4.2}
Under the conditions and notation of Theorem~\ref{Th_oiv_4.1}:
\begin{itemize}
\item[1)] for ${p>0}$, ${q>0}$, in some neighborhood of
$y_0(\cdot)$ in $C^1[a;b]$ the estimate
\begin{equation}\label{oiv_19}
\Phi(y)-\Phi(y_0)\geq \frac{1}{2}\min(p,q)\cdot \|y\|^2_{H^1[a;b]}
\end{equation}
holds;

\item[2)] for ${p>0}$, ${q\geq 0}$, in some neighborhood of
$y_0(\cdot)$ in $C^1[a;b]$ the estimate
\begin{equation}\label{oiv_20}
\Phi(y)-\Phi(y_0)\geq \frac{\pi^2p}{2(\pi^2+16(b-a)^2)}\cdot
\|y\|^2_{H^1[a;b]}
\end{equation}
holds;

\item[3)] for ${p>0}$, ${q<0}$, under the
restriction~\eqref{oiv_18} to a length of ${[a;b]}$, in some
neighborhood of $y_0(\cdot)$ in $C^1[a;b]$ the estimate
\begin{equation}\label{oiv_21}
\Phi(y)-\Phi(y_0)\geq
\frac{\pi^2p-16(b-a)^2|q|}{2(\pi^2+16(b-a)^2)}\cdot
\|y\|^2_{H^1[a;b]}
\end{equation}
holds.
\end{itemize}
\end{theorem}

At last, applying Theorem~\ref{Th_oiv_3.1} to the auxiliary
integrand $\widetilde{f}$ leads to solution of the inverse extreme
problem for $\Phi$ at an arbitrary point ${y_0(\cdot)\in
C^2[a;b]}$.

\begin{theorem}\label{Th_oiv_4.3}
Let, under the conditions of Theorem~\ref{Th_oiv_4.1}, the
variational functional~\eqref{oiv_1} attains a local minimum at a
point ${y_0(\cdot)\in C^2[a;b]}$ under the boundary conditions
${y(a)=y_0(a)}$, ${y(b)=y_0(b)}$ and the strengthened Legendre
condition. Then the integrand $f$ takes form of
$$
f(x,y,z)=P(x,y-y_0(x))+\biggl(C+\int\limits_a^x
P_y(t,-y_0(t))dt+[q(x,y-y_0(x))-q(x,-y_0(x))]\biggr)\cdot
$$
\begin{equation}\label{oiv_22}
\cdot(z-y'_0(x))+\frac{1}{2}\biggl(p(x)+[\rho(x,y-y_0(x),z-y'_0(x))-\rho(x,-y_0(x),-y'_0(x))]\biggr)\cdot
(z-y'_0(x))^2\, ,
\end{equation}
where ${C\in\mathbb{R}}$; $P$, $q$, ${p>0}$, ${\rho\in C^2}$  can
be chosen arbitrarily.
\end{theorem}

From here a formula of the general form of the
functional~\eqref{oiv_1} taking a local minimum in $C^1[a;b]$ at a
point ${y_0(\cdot)\in C^2[a;b]}$ under the strengthened Legendre
condition:
$$
\Phi(y)=\int\limits_a^b
\biggl(P(x,y-y_0(x))+\biggl[\int\limits_a^x
P_y(t,-y_0(t))dt+q(x,y-y_0(x))-q(x,-y_0(x))\biggr]\cdot
$$
\begin{equation}\label{oiv_23}
\cdot(y'-y'_0(x))+\frac{1}{2}\biggl[p(x)+\rho(x,y-y_0(x),y'-y'_0(x))-\rho(x,-y_0(x),-y'_0(x))\biggr]\cdot
(y'-y'_0(x))^2\biggr)dx\, ,
\end{equation}
where ${C\in\mathbb{R}}$; $P$, $q$, ${p>0}$, ${\rho\in C^2}$  can
be chosen arbitrarily, arises.

Thus, under the strengthened Legendre condition, the inverse
extreme variational problem at an arbitrary point ${y_0(\cdot)\in
C^2}$ is solved: the all functionals of the type~\eqref{oiv_1},
attaining a local minimum at a point ${y_0(\cdot)}$, are
described.

\section{\textbf{Case of compact extremum in $H^1[a;b]$}}

In the Hilbert--Sobolev space ${W^{1,2}[a;b]=H^1[a;b]}$ equipped
with the norm
\begin{equation}\label{oiv_24}
\|y\|^2_{H^1[a;b]}=\int\limits_a^b (y^2+y'^2)dx\,\,,
\end{equation}
as it is well known, by virtue of I.V. Skrypnik
theorem~(\cite{oiv_4}, Ch.11) the nonabsolute local extrema of the
variational functionals practically absend. Note that in the
present work the norm~\eqref{oiv_24} was appeared above
(Theorem~\ref{Th_oiv_2.1},~\ref{Th_oiv_4.2}) by natural way even
for extreme problems in $C^1[a;b]$.

In the our works~\cite{oiv_5}--\cite{oiv_7} and in the works by
E.V. Bozhonok~\cite{oiv_8}--\cite{oiv_10} a general notion of
\textit{compact extremum} (or \textit{$K$--extremum}) of a
functional was studied (see, also,~\cite{oiv_11}). It have been
shown there that the classical, both necessary and sufficient
conditions of a local extremum of variational functional in
$C^1[a;b]$ can be extended to the case of $K$--extremum in
$H^1[a;b]$. In this case, $K$--extrema inherit the important
properties of the local extrema and can be considered as an analog
of the ones in the case of variational functionals in $H^1[a;b]$.
Let's bring a relevant information.

\begin{definition}
Let a real functional ${\Phi:H\rightarrow \mathbb{R}}$ be defined
in a Hilbert space $H$. Say that $\Phi$ has a \textit{compact
minimum} (or \textit{$K$--extremum}) at a point ${y_0\in H}$ if,
for each absolutely convex (a.c.) compact set ${C\subset H}$, the
restriction of $f$ to the subspace ${(y_0+span\,C)}$ has a local
minimum at $y_0$ respective to Banach norm ${\|\cdot\|_C}$ in
$span\,C$ generated by $C$. In other words, for each a.c.
compactum ${C\subset H}$ there exists such
${\varepsilon=\varepsilon(C)>0}$ that
${\varphi(y)\geq\varphi(y_0)}$ as ${y-y_0\in \varepsilon\cdot C}$.
\end{definition}

The well posedness and the validity for the case of $K$--extremum
of the variational functional~\eqref{oiv_1} of the classical
extreme conditions in $C^1$ (Euler--Lagrange equation, Legendre
condition, Jacobi condition) require, as it was shown
in~\cite{oiv_7}, belonging coefficient $R(x,y,z)$ in the quadratic
representation~\eqref{oiv_10} of the integrand $f$:
$$
f(x,y,z)=P(x,y)+Q(x,y)\cdot z+\frac{1}{2}R(x,y,z)\cdot z^2
$$
to an appropriate \textit{dominated mixed smoothness space}
$C^2_{xy}$ (see~\cite{oiv_12},~\cite{oiv_13}). Namely, for the
arbitrary compacta ${C_x,\,\,C_y\subset\mathbb{R}}$ the following
property holds:
$$
(x\in C_x\,,\,\,y\in C_y\,,\,\,-\infty<z<+\infty)\Rightarrow
(R(x,y,z)\,\, is \,\, uniformly\,\, continuous \,\, and
$$
$$
bounded,\,\,together\,\, with\,\, its\,\, first\,\, and\,\,
second\,\, partial\,\, derivatives).
$$

Under the conditions above, the Euler--Lagrange equation, Legendre
condition, strengthened Legendre condition and Jacobi condition
for the Hamilton--Jacobi equation are extended to the case of
$K$--extremum in an arbitrary $W^{2,2}$--smooth point
${y_0(\cdot)\in H^1[a;b]}$. It allows to extend the results of i.4
to the case of $K$--minimum in ${H^1[a;b]}$. Let's bring the
corresponding formulations.

\begin{theorem}\label{Th_oiv_5.1}
Let the variational functional~\eqref{oiv_1} at a
$W^{2,2}$--smooth point ${y_0(\cdot)\in H^1[a;b]}$ satisfies
Euler--Lagrange equation~\eqref{oiv_17}, in addition ${R(x,y,z)\in
C^2_{xy}}$. Then, under the conditions and notation of
Theorem~\ref{Th_oiv_4.1},
\begin{itemize}
\item[1)] for ${p>0}$, ${q\geq 0}$, $\Phi(y)$ attains a strong
$K$--minimum at $y_0(\cdot)$ (without any restriction to a length
of ${[a;b]}$);

\item[2)]for ${p>0}$, ${q<0}$, and under the
restriction~\eqref{oiv_18} to a length of ${[a;b]}$, $\Phi(y)$
attains a strong $K$--minimum at $y_0(\cdot)$ as well.
\end{itemize}
\end{theorem}

\begin{theorem}\label{Th_oiv_5.2}
Under the conditions and notation of Theorem~\ref{Th_oiv_5.1}:
\begin{itemize}
\item[1)] for ${p>0}$, ${q>0}$, for each a.c. compactum ${C\subset
H^1[a;b]}$ there exists such ${\varepsilon=\varepsilon(C)>0}$ that
inclusion ${y-y_0\in \varepsilon\cdot C}$ implies
 estimate~\eqref{oiv_19};

\item[2)] for ${p>0}$, ${q\geq 0}$, for each a.c. compactum
${C\subset H^1[a;b]}$ there exists such
${\varepsilon=\varepsilon(C)>0}$ that inclusion ${y-y_0\in
\varepsilon\cdot C}$ implies
 estimate~\eqref{oiv_20};

\item[3)] for ${p>0}$, ${q<0}$, under the
restriction~\eqref{oiv_18} to a length of ${[a;b]}$, for each a.c.
compactum ${C\subset H^1[a;b]}$ there exists such
${\varepsilon=\varepsilon(C)>0}$ that inclusion ${y-y_0\in
\varepsilon\cdot C}$ implies
 estimate~\eqref{oiv_21}.
\end{itemize}
\end{theorem}

\begin{theorem}\label{Th_oiv_5.3}
Let, under the conditions and notation of
Theorem~\ref{Th_oiv_5.1}, the variational functional~\eqref{oiv_1}
attains a $K$--minimum at a $W^{2,2}$--smooth point ${y_0(\cdot)}$
from ${H^1[a;b]}$ under the boundary conditions ${y(a)=y_0(a)}$,
${y(b)=y_0(b)}$ and under the strengthened Legendre condition.
Then the integrand $f$ takes form of~\eqref{oiv_22}, where
${C\in\mathbb{R}}$; $P$, $q$, ${p>0}$ from $C^2$  and $\rho$ from
$C^2_{xy}$ can be chosen arbitrarily.
\end{theorem}

From here the formula~\eqref{oiv_23} of the general form of the
functional~\eqref{oiv_1} having a $K$--minimum at a
$W^{2,2}$--smooth point ${y_0(\cdot)}$ from ${H^1[a;b]}$ under the
strengthened Legendre condition, follows.

\end{document}